\documentclass[final,leqno]{siamltex}

\usepackage{amsfonts}
\usepackage{epsfig}
\usepackage{amsmath}
\title{Delay-induced patterns in a two-dimensional lattice of coupled oscillators}

\author{
Markus Kantner\thanks{Weierstrass Institute for Applied Analysis and Stochastics, Mohrenstr. 39, 10117 Berlin, Germany ({\tt kantner@wias-berlin.de})},
\and Serhiy Yanchuk\thanks{Humboldt University of Berlin, Institute of Mathematics, Unter den Linden 6, 10099 Berlin, Germany ({\tt yanchuk@math.hu-berlin.de})},
\and Eckehard Sch{\"o}ll\thanks{Technical University of Berlin, Institute of Theoretical Physics, Hardenbergstr. 36, 10623 Berlin, Germany ({\tt schoell@physik.tu-berlin.de}).}
}

\begin{document}

\maketitle

\begin{abstract}
We show how a variety of stable spatio-temporal periodic patterns
can be created in 2D-lattices of coupled oscillators with non-homogeneous
coupling delays. A ``hybrid dispersion relation'' is introduced,
which allows studying the stability of time-periodic patterns analytically
in the limit of large delay. The results are illustrated using the
FitzHugh-Nagumo coupled neurons as well as coupled limit cycle (Stuart-Landau)
oscillators.
\end{abstract}

\begin{keywords}
time delay, patterns , FitzHugh-Nagumo neurons,  lattices, componentwise
timeshift-transformation
\end{keywords}

\begin{AMS}
34K13, 34K18, 34K20, 34K35, 34K60, 37L60, 37N25, 37N35, 93C23
\end{AMS}

\pagestyle{myheadings}
\thispagestyle{plain}
\markboth{M. KANTNER, S. YANCHUK AND E. SCH\"OLL}{DELAY-INDUCED PATTERNS IN A 2D-LATTICE}

\section{Introduction\label{sec: model system}}

Coupled dynamical systems with time delays arise in various applications
including semiconductor lasers \cite{Franz2008,Ludge2012,Soriano2013,Williams2013b}, electronic circuits \cite{Reddy2000}, optoelectronic oscillators \cite{Williams2013}, neuronal networks \cite{Izhikevich2006,Stepan2009,Campbell2007}, gene regulation networks \cite{Tiana2013}, socio-economic systems \cite{Szalai2013,Orosz2010}
and many others \cite{Fiedler2007,Flunkert2013,Zakharova2013,Kinzel2009,Flunkert2010b,Just2010}.
Understanding the dynamics in such systems is a challenging task.
Even a single oscillator with time delayed feedback exhibits phenomena,
which are not expected in this class of systems, such as Eckhaus instability
\cite{Wolfrum2006}, coarsening \cite{Giacomelli2012}, or chimera
state \cite{Larger2013}. Some of them, like low frequency fluctuations
in laser systems with optical feedback are still to be understood
\cite{Pieroux2003}. The situation is even more complicated when
several systems are interacting with non-identical delays. In this
case, somewhat more is known about some specific coupling configurations,
e.g. ring \cite{Popovych2011,Yanchuk2011,Choe2014,Bungay2007}, and less on more complex
coupling schemes \cite{Izhikevich2006,Luecken2013b,DHuys2013,Dahms2012,Cakan2014}. Recently,
it has been shown that a ring of delay coupled systems possesses a
rich variety of stable spatiotemporal patterns \cite{Popovych2011,Yanchuk2011}.
For the neuronal models in particular, this implies the existence
of a variety of spiking patterns induced by the delayed synaptic connections.

Here we present a system with time delayed couplings, which is capable
of producing a variety of stable two-dimensional spatio-temporal patterns.
More specifically, we show that a 2D regular set of dynamical systems
$\mathbf{u}_{m,n}(t)$ (neuronal models can be used) may exhibit a stable periodic behavior (periodic spiking) such that the oscillator $\mathbf{u}_{m,n}(t)$
reaches its maximum (spikes) at a time $\eta_{m,n}$, which can be
practically arbitrary chosen within the period. For this, time delays
should be selected accordingly to some given simple rule. As particular
cases, the synchronous, cluster, or splay states can be realized. 

Our work is a generalization of the previous results on the ring \cite{Popovych2011,Yanchuk2011},
extending them to the two-dimensional case. However, the analysis,
which we have to employ has important differences. In particular,
the combination of the spatial structure of the system (spatial coordinates
$m$ and $n$) and temporal delays required the introduction of a
so called ``hybrid dispersion relation'' for the investigation of
the stability of stationary state and nonlinear plane waves in the
homogeneous system. Roughly speaking, this hybrid dispersion relation
is a synthesis of the dispersion relation from the pattern formation
theory in spatially extended systems \cite{Cross1993,Cross2009}
and the pseudo-continuous spectrum developed for purely temporal delay
systems \cite{Wolfrum2006,Lichtner2011,Wolfrum2010}. 

We believe that such a higher-dimensional extension allows to speculate
about the possibility of employing such systems for generating or
saving visual patterns, and can be probably of use for information
processing purposes.
Small arrays of delay-coupled optoelectronic oscillators have indeed already been realized experimentally \cite{Williams2013}. Similarly, autonomous Boolean networks of electronic logic gates have been demonstrated as versatile tools for the realizations of various space-time patterns \cite{Rosin2013}.
Moreover, our analysis provides another evidence
that the delays in coupled systems can play a constructive functional
role.

\begin{figure}[h]
\centering
\includegraphics{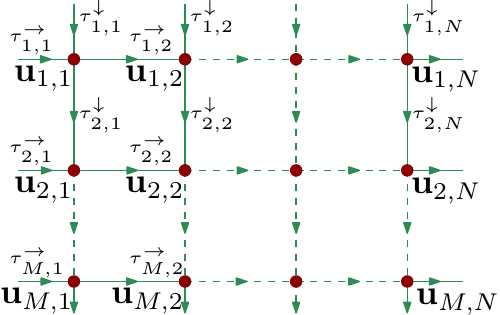} \caption{Coupling scheme. The dynamics of each node $\mathbf{u}_{m,n}(t)$
is described by system (\ref{eq: general system}). Each coupling
connection possesses a delay $\tau_{m,n}^{\downarrow}$ or $\tau_{m,n}^{\rightarrow}$.
All edges are unidirectional.\label{pic: 2-dim lattice coupling structure}}
\end{figure}

More specifically, we consider a lattice of $M\times N$ delay-coupled
systems (delay differential equations) of the form 
\begin{align}
\frac{\mathrm{d}}{\mathrm{d}t} \mathbf{u}_{m,n}(t) &= \mathbf{F}\left(\mathbf{u}_{m,n}(t),\mathbf{u}_{m-1,n}(t-\tau_{m,n}^{\downarrow})+\mathbf{u}_{m,n-1}(t-\tau_{m,n}^{\rightarrow})\right),\label{eq: general system}\\
 &\qquad m=1,\ldots ,M; \quad n=1,\dots, N, \nonumber 
\end{align}
where $\mathbf{F}:\mathbb{R}^{d}\times\mathbb{R}^{d}\to\mathbb{R}^{d}$
is a nonlinear function determining the dynamics of $\mathbf{u}_{m,n}\in\mathbb{R}^{d}$
in the lattice. The indices $m$ and $n$ determine the position of
the node, see Fig. \ref{pic: 2-dim lattice coupling structure}. We
assume periodic boundary conditions $\mathbf{u}_{M+1,n}\equiv\mathbf{u}_{1,n}$
and $\mathbf{u}_{m,N+1}\equiv\mathbf{u}_{m,1}$ such that the system
has translation invariance. Time delays $\tau_{m,n}^{\downarrow}$
and $\tau_{m,n}^{\rightarrow}$ denote the connection delays between
the corresponding nodes. Since each node has two incoming connections,
the arrows $\downarrow$ and $\rightarrow$ correspond to the coupling
from the node located above, respectively left, see Fig.~\ref{pic: 2-dim lattice coupling structure}.
Here we restrict the analysis to two systems: Stuart-Landau (SL) oscillators
as a simplest dynamical system exhibiting limit cycle behavior and
FitzHugh-Nagumo (FHN) systems as a representative of conductance based,
excitable neuronal models \cite{FitzHugh1961,Nagumo1962}. While the first
model allows for a deeper analytical insight, the second one can be
studied mainly numerically and shows similar qualitative results. 

\begin{figure}
\includegraphics[width=1\textwidth]{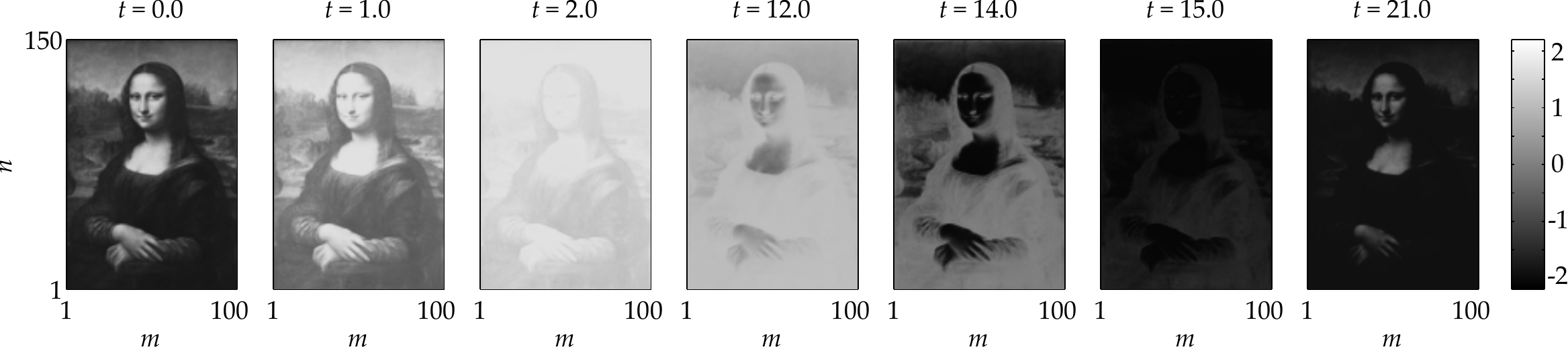}\caption{\label{fig:ML}Snapshots of the spatio-temporal behavior in a system
of $100\times150$ identical FHN neurons Eq. (\ref{eq:FHN}) with appropriately
adjusted time delays $\tau_{m,n}^{\downarrow}$ and $\tau_{m,n}^{\rightarrow}$.
At each grid point with the coordinate $(m,n)$, the level of gray
(see colorbar) corresponds to the membrane voltage $v_{m,n}(t)$ at
this time moment. The pattern reappears periodically with a time period
$T=21.95$. More details are given in Sec.~\ref{sec:non}. }
\end{figure}

An example of a stable spatio-temporal pattern in a lattice of $100\times150$
FHN neurons with non-homogeneous delays, the ``Mona Lisa''-pattern
is shown in Figure~\ref{fig:ML}. Each frame corresponds to a snapshot
at a fixed time $t$ and the different level of gray at a point $(m,n)$
corresponds to the value of the voltage component of $\mathbf{u}_{m,n}(t)$
at this time $t$. More details on how such patterns can be created
are given in the following sections, especially in Sec.~\ref{sec:non}.

The structure of the remaining part of the paper is as follows: In
Sec.~\ref{sec:SLhom} we consider SL systems with homogeneous time
delays $\tau_{m,n}^{\rightarrow}=\tau_{m,n}^{\downarrow}=\tau$. We
investigate the stability of the homogeneous steady state as well
as various plane wave solutions in the system. The number of stable
plane wave solutions is shown to increase with the delay. Further
in Sec.~\ref{sec:FHN-hom} similar results are obtained for the FHN
systems. Section \ref{sec:non} considers the case when the delays
are not identical. In this case it is shown how a variety of spatiotemporal
patterns can be created by varying the coupling delays. Finally, additional
illustrative examples are presented in Sec.~\ref{sub:Ex}.

\section{Stuart-Landau oscillators with homogeneous coupling delays\label{sec:SLhom}}

In this section we start with a lattice of SL oscillators with homogeneous
delays $\tau_{m,n}^{\rightarrow}=\tau_{m,n}^{\downarrow}=\tau$ as
a model for the local dynamics: 
\begin{align}
\begin{aligned}
\frac{\mathrm{d}}{\mathrm{d}t}z_{m,n}(t) &= (\alpha+i\beta)z_{m,n}(t) - z_{m,n}(t)|z_{m,n}(t)|^{2}+ \\
&\quad + \frac{C}{2}\left(z_{m-1,n}(t-\tau)+z_{m,n-1}(t-\tau)\right).\label{eq: SL}
\end{aligned}
\end{align}
The variables $z_{m,n}$ are complex-valued. The parameter $\alpha$ controls
the local dynamics without coupling, i.e. a stable steady state exists
for $\alpha<0$ and a stable limit cycle for $\alpha>0$; $\beta$
is the frequency of this limit cycle. The coupling strength is determined
by $C>0$. 

It this section we study the bifurcation scenario, which is associated
with the destabilization of the homogeneous steady state $z=0$ and
the appearance of various plane waves. Many aspects of this scenario
can be studied analytically due to the $S^{1}$ equivariance of the
system: $\mathbf{F}(e^{i\nu}\cdot,e^{i\nu}\cdot)=e^{i\nu}\mathbf{F}(\cdot,\cdot)$
for any real $\nu$. At some places we will assume additionally that
the delay $\tau$ is large comparing to the timescale of the system
(we will mention it each time explicitly), which simplifies analytical
calculations.

\subsection{Stability and bifurcations of homogeneous stationary state}

\subsubsection{Characteristic equation and its solutions}

System (\ref{eq: SL}) has a homogeneous steady state $z_{m,n}\equiv0$.
Let us investigate stability of this stationary state. Linearizing
the equation of motion (\ref{eq: SL}) around $z_{m,n}=0$ yields
the following equation for the evolution of small perturbations $\delta z_{m,n}(t)$:
\[
\frac{\mathrm{d}}{\mathrm{d}t}\delta z_{m,n}(t)=(\alpha+i\beta)\delta z_{m,n}(t)+\frac{C}{2}\left(\delta z_{m-1,n}(t-\tau)+\delta z_{m,n-1}(t-\tau)\right).
\]
This equation can be diagonalized by a spatial discrete Fourier-transformation
\[
\delta\tilde{z}_{k_{1},k_{2}}=\sum_{m=1}^{M}\sum_{n=1}^{N}e^{ik_{1}m+ik_{2}n}\delta z_{m,n},
\]
where the wavevector $\mathbf{k=}(k_{1},k_{2})^{T}$ admits the following
discrete values: 
\begin{equation}
(k_{1},k_{2})^{T}=2\pi\left(\frac{l}{M},\frac{j}{N}\right)^{T},\quad l=1,\ldots,M,\quad j=1,\ldots,N.\label{eq:k}
\end{equation}
We obtain
\begin{equation}
\frac{\mathrm{d}}{\mathrm{d}t}\delta\tilde{z}_{k_{1},k_{2}}(t)=(\alpha+i\beta)\delta\tilde{z}_{k_{1},k_{2}}(t)+\frac{C}{2}\left(e^{ik_{1}}+e^{ik_{2}}\right)\delta\tilde{z}_{k_{1},k_{2}}(t-\tau).\label{eq:SLld}
\end{equation}
Since the equation for the Fourier modes is uncoupled, one can drop
the indices $k_{1}$ and $k_{2}$ for simplicity and introduce the
notations 
\begin{equation}
k_{\pm}=\frac{1}{2}\left(k_{1}\pm k_{2}\right),\label{eq:Q}
\end{equation}
which is just a rotation of coordinates in the Fourier space. Note
that $k_{\pm}$ admits discrete values accordingly to (\ref{eq:k}).
Then the system (\ref{eq:SLld}) can be rewritten as 
\[
\frac{\mathrm{d}}{\mathrm{d}t}\delta\tilde{z}(t)=(\alpha+i\beta)\delta\tilde{z}(t)+Ce^{ik_{+}}\cos\left(k_{-}\right)\delta\tilde{z}(t-\tau).
\]
Since the modes are decoupled in the Fourier space, the corresponding
characteristic equation factorizes and reads

\begin{equation}
0=\prod_{\{(k_{_{+}},k_{-})\}} \left( -\lambda+\alpha+i\beta+C\cos\left(k_{-}\right)e^{ik_{+}-\lambda\tau} \right).\label{eq: SLststChEq-1}
\end{equation}
This is a transcendental equation in $\lambda$, which can be solved
as 
\begin{equation}
\lambda_{j}=\alpha+i\beta+\frac{1}{\tau}W_{j}\left(\tau C\cos\left(k_{-}\right)e^{ik_{+}-(\alpha+i\beta)\tau}\right),\label{eq:W}
\end{equation}
for all admitted pairs ${(k_{+},k_{-})}$, where $j$ denotes the
$j$th branch of the Lambert $W$-function. This formula can be used
for studying the stability of the homogeneous state. In particular,
if all eigenvalues $\lambda_{j}$ have negative real parts for all
possible wavevectors $k_{\pm}$, then the steady state is asymptotically
stable.

\subsubsection{Large delay approximation\label{sub:ststLong}}

\begin{figure}
\centering\includegraphics[width=0.34\textwidth]{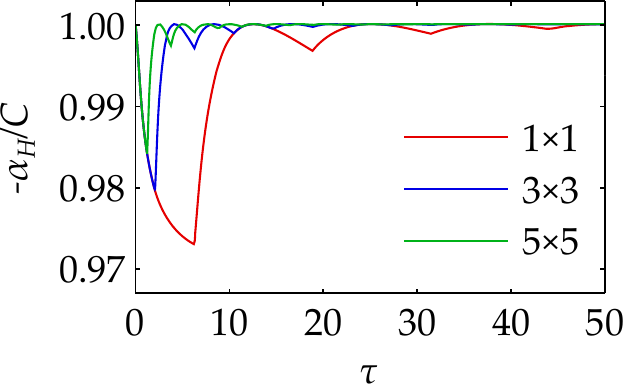}

\caption{Destabilization threshold of the homogeneous stationary state in lattices
of coupled Stuart-Landau oscillators as a function of the coupling delay. The
asymptotic result $\alpha_{H}=-C$ is exact for vanishing and infinite
delay. However, also for moderate delays the approximation sufficiently
matches the numerical results. \label{pic:SLfirstHopf}}
\end{figure}

A deeper understanding of the spectrum can be achieved for large delays
using the asymptotic methods \cite{Lichtner2011,Sieber2013,Flunkert2010b}. Accordingly
to these results, the spectrum splits generically into two parts accordingly
to the asymptotics for large delays. The first part is called the
\emph{strongly unstable spectrum} and the second part is the \emph{pseudo-continuous
spectrum}. The strong spectrum exists for $\alpha>0$ and consists
of two complex conjugate, isolated roots which are close to $\lambda_{S,\pm}=\alpha\pm i\beta$.
In such a case, the contribution of the term $\frac{1}{\tau}W_{j}\left[\cdot\right]$
in (\ref{eq:W}) vanishes. In fact, the strong spectrum always converges
to the unstable part of the spectrum of the system with omitted delayed
terms \cite{Yanchuk2005,Wolfrum2010,Lichtner2011}, i.e. $\frac{\mathrm{d}}{\mathrm{d}t}\delta\tilde{z}(t)=(\alpha+i\beta)\delta\tilde{z}(t)$
in this case. Besides the strong spectrum there are infinitely many
more eigenvalues, accumulating on some curves in the complex plane
as $\tau\to\infty$. These eigenvalues form the pseudo-continuous
spectrum and can be found by substituting the ansatz $\lambda=\frac{1}{\tau}\gamma(\Omega)+i\Omega$
in the characteristic equation (\ref{eq: SLststChEq-1}). One obtains
\[
i\Omega=\alpha+i\beta+C\cos\left(k_{-}\right)e^{ik_{+}}e^{-\gamma-i\Omega\tau},
\]
where the small term $\gamma/\tau$ has been neglected. It can be
solved as 
\[
Y:=e^{-\gamma}e^{i(k_{+}-\Omega\tau)}=-\frac{\alpha+i(\beta-\Omega)}{C\cos{\left(k_{-}\right)} }
\]
and finally for the curve $\gamma(\Omega)$ we obtain 
\begin{equation}
\gamma(\Omega)=-\frac{1}{2}\log|Y|^{2}=-\frac{1}{2}\log\left(\frac{\alpha^{2}+(\beta-\Omega)^{2}}{C^{2}\cos^2{\left(k_{-}\right)} }\right).\label{eq: SLststpcs}
\end{equation}
Note that the spectrum is invariant with respect to complex conjugation.
It is easy to see that the spatial mode $k_{-}=0$ corresponds to
the maximal values of $\gamma(\Omega)$. Thus, the spatial modes with
$k_{-}=0$ are most unstable. Moreover, it is easy to check that the
leading value of $\gamma$ is negative if $\left| \alpha \right|>C$
and positive otherwise. This implies that the homogeneous steady state
is asymptotically stable for $\alpha<-C$ and unstable for $\alpha>-C$
(we take also into account that there is a strongly unstable spectrum
for $\alpha>0$). Hence, the destabilization takes place at $\alpha=-C$
via Hopf bifurcation at the frequency $\Omega\approx\beta$. Here
we should remark that our conclusions are rigorous for long enough
delay $\tau$. In our case, however, the destabilization threshold
is well approximated by this value already for moderate delays, see
Fig. \ref{pic:SLfirstHopf}.

\begin{figure}[h]
%\centering \includegraphics[height=2.4cm]{fig/SL_pcs_stst_destabilization_new2}\includegraphics[height=2.4cm]{fig/SL_2torus_stst_HDR}
\includegraphics[width=\textwidth]{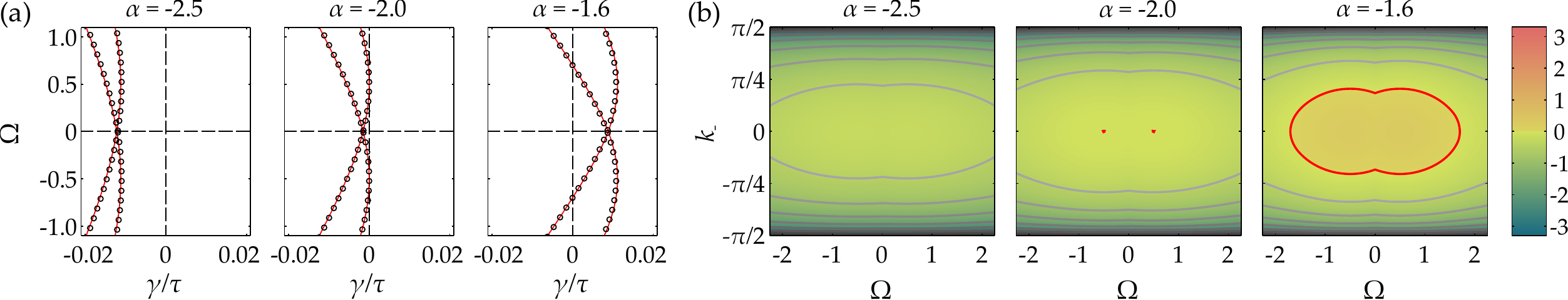}
\caption{(a) Destabilization of the zero steady state in $3\times3$ lattice
of delay coupled SL oscillators with $C=2$, $\beta=0.5$, and $\tau=20$.
The plot shows numerically computed eigenvalues (solutions of Eq.
(\ref{eq: SLststChEq-1})) and the continuous large delay approximation
(\ref{eq: SLststpcs}) by the red line. At $\alpha=-2.5$ the stationary
state is stable, $\alpha=-2$ is the critical case, and for $\alpha=-1.6$
an unstable steady state is shown. All eigenvalues accumulate along
the curves $\gamma(\Omega,k_{-})$ given by Eq. (\ref{eq: SLststpcs})
with maxima at $\Omega=\pm\beta$. (b) Destabilization of the stationary
state in an infinitely large lattice, where $k_{-}$becomes a continuous
parameter. The density plot shows the scaled real part of the eigenvalues $\operatorname{Re} \gamma = \lambda/\tau$. The parameters are the same as on the left. \label{pic: SL-eig-stst}}
\end{figure}

The illustration by the numerically computed eigenvalues is shown
in Fig.~\ref{pic: SL-eig-stst} for the system of $3\times3$ coupled
SL oscillators for three cases: stable, critical, and unstable. One
can observe also how multiple Hopf bifurcation may emerge after the
destabilization. In the following section we discuss the plane waves
arising in these Hopf bifurcations.

\subsection{Nonlinear plane waves and their stability}

\subsubsection{Branches of plane wave solutions}

\begin{figure}
\centering
\includegraphics[scale=0.75]{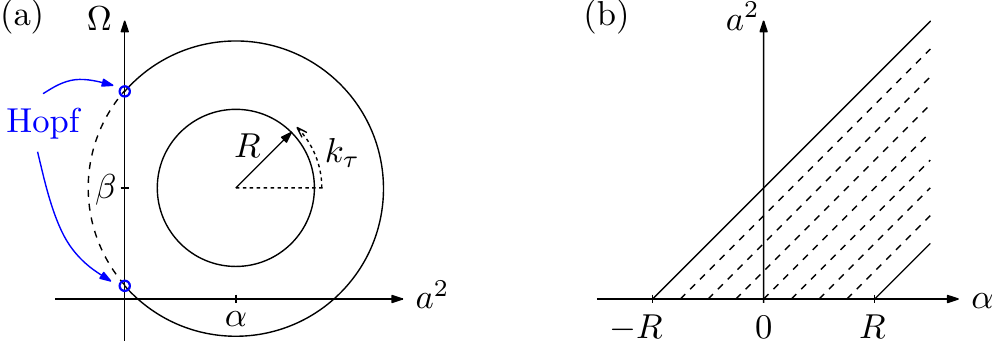}
\caption{(a) Frequency $\Omega$ versus amplitude $a$ of periodic solution
(\ref{eq:persol}). The admissible parameters are located on the circles
with radii $R$ {[}Eq.~\ref{pic: SLpsolEll}{]}, except the parts
with $a^{2}<0$. (b) Branches of periodic solutions emerge from the
steady state within the interval $\alpha\in\left[-\left|R\right|,\left|R\right|\right]$.
\label{pic: SLpsolEll}}
\end{figure}

Because of the $S^{1}$-equivariance of the Stuart-Landau system (\ref{eq: SL}),
periodic solutions emerging from the homogeneous steady state via
Hopf bifurcations have the following explicit form 
\begin{equation}
z(t)=ae^{i\Omega t-ik_{1}m-ik_{2}n}.\label{eq:persol}
\end{equation}
The amplitude $a$, frequency $\Omega$, and the wavevector $\mathbf{k}=(k_{1},k_{2})^{T}$
of the periodic solutions satisfy the following relation 
\[
i\Omega=\alpha+i\beta-a^{2}+e^{i(k_{+}-\Omega\tau)}C\cos{\left(k_{-}\right)},
\]
which is obtained by substituting (\ref{eq:persol}) into (\ref{eq: SL}).
Here we recall that $k_{\pm}$ is the wavevector in rotated coordinates
(\ref{eq:Q}). Taking real and imaginary parts, we obtain 
\begin{subequations}
\begin{align}
a^{2} & =\alpha+R\cos k_{\tau},\label{eq: SLpsolRe}\\
\Omega & =\beta+R\sin k_{\tau},\label{eq: SLpsolIm}
\end{align}
\end{subequations}
where we denote $R:=C\cos k_{-}$ and $k_{\tau}:=k_{+}-\Omega\tau$.
By excluding $k_{\tau}$ we obtain 
\begin{equation}
(a^{2}-\alpha)^{2}+(\Omega-\beta)^{2}=R^{2}.\label{eq: SLpsolCir}
\end{equation}
Therefore all periodic solutions can be found on circles in the $(a^{2},\Omega)$-parameter
space as shown in Fig. \ref{pic: SLpsolEll}. Equation (\ref{eq: SLpsolIm})
is known as \emph{Kepler's equation} and can be solved numerically
with respect to $\Omega$. The number of solutions of (\ref{eq: SLpsolIm})
matches the number of Hopf bifurcations and periodic solutions. All
possible frequencies are confined to the interval $-\left|R\right|+\beta\leq\Omega\leq\left|R\right|+\beta$.
The corresponding value for the parameter $\alpha$, at which Hopf
bifurcation occurs, is then given by Eq. (\ref{eq: SLpsolRe}) with
$a^{2}=0$. For a system of the size $N\times N$, the number of Hopf
bifurcations, or periodic solutions respectively, can be estimated
by 
\begin{eqnarray*}
\mathcal{N}_{H} \approx \begin{cases}
\frac{2C\tau}{\pi} N \frac{1}{\sin{\left(\frac{\pi}{2 N}\right)}} & N\,\mathrm{even},\\
\frac{2C\tau}{\pi} N \cot{\left( \frac{\pi}{2N} \right)} & N\,\mathrm{odd}.
\end{cases}
\end{eqnarray*}
In the case of large $N$ the asymptotic number of solutions is given
by 
\begin{equation}
\mathcal{N}_{H}\sim\frac{4C\tau}{\pi^{2}}N^{2}.\label{eq:Npsol}
\end{equation}
Thus, in the case of large delay or lattice-size the number of solutions
grows and any point on the circles (\ref{eq: SLpsolCir}) refers to
a periodic solution, i.e. the circle disk is densely filled with points
$(a,\Omega)$ corresponding to the existing periodic solutions. Indeed,
for an infinitely large lattice, the wavevector $\mathbf{k}$ admits
any values from $\mathbb{R}^{2}$. Therefore one obtains continuous
parameters $k_{-}\in\left(-\pi,\pi\right)$ and $k_{\tau}=(k_{+}-\Omega\tau)\in\left(0,2\pi\right)$,
which parametrize all periodic solutions of Eq. (\ref{eq: SL}). 

\subsubsection{Stability of plane wave solutions}

The local asymptotic stability of plane wave solutions can be studied
using the linearized equation for small perturbations $\delta z_{m,n}(t)$.
In co-rotating coordinates $\xi_{m,n}$ determined by 
\begin{equation}
z_{m,n}(t)=e^{i(\Omega t-k_{1}m-k_{2}n)}u_{m,n}(t),\label{eq:COR}
\end{equation}
the plane wave is just the steady state $u_{m,n}=a$. Therefore, after
changing the coordinates (\ref{eq:COR}), substituting $u_{m,n}(t)=a+\xi_{m,n}(t)$,
and linearizing the obtained equation in small perturbations $\xi_{m,n}(t)$,
we obtain 
\begin{align}
\begin{aligned}  \label{eq:xi}
\frac{\mathrm{d}}{\mathrm{d}t}\xi_{m,n}(t)  &= \left(\alpha+i\beta-i\Omega-2a^{2}\right)\xi_{m,n}(t)-a^{2}\xi_{m,n}^{\ast}(t) + \\
 &\quad +\frac{C}{2}e^{-i\Omega\tau}\left(e^{ik_{1}}\xi_{m-1,n}(t-\tau)+e^{ik_{2}}\xi_{m,n-1}(t-\tau)\right).
\end{aligned}
\end{align}
The solutions of the obtained linear equation can be found by the
ansatz
\begin{equation}
\xi_{m,n}(t)=b_{1}e^{\lambda t-iq_{1}m-iq_{2}n}+b_{2}^{\ast}e^{\lambda^{\ast}t+iq_{1}m+iq_{2}n}.\label{eq:bloch}
\end{equation}
The ansatz (\ref{eq:bloch}) can be obtained, e.g. by rewriting the
system (\ref{eq:xi}) in the real form, diagonalizing it with the
discrete Fourier transform (similar to (\ref{eq:SLld})), and noticing
that the equations for the Fourier components $\tilde{\xi}_{q_{1},q_{2}}$
and $\tilde{\xi}_{-q_{1},-q_{2}}$ are complex conjugate, and hence,
they have the same stability properties with the complex conjugate
eigenvalues, see also \cite{Cross2009}. After substituting (\ref{eq:bloch})
into (\ref{eq:xi}), the coefficients at the two linearly independent
functions $e^{\lambda t-iq_{1}m-iq_{2}n}$ and $e^{\lambda^{*}t+iq_{1}m+iq_{2}n}$
should be zero. This leads to the system of two linear equations for
unknowns $b_{1}$ and $b_{2}$: 
\begin{align*}
\left[\alpha+i\beta-i\Omega-\lambda-2a^{2}+\frac{C}{2}e^{-(\lambda+i\Omega)\tau}\left(e^{i(k_{1}+q_{1})}+e^{i(k_{2}+q_{2})}\right)\right]b_{1}-a^{2}b_{2} & =0,\\
-a^{2}b_{1}+\left[\alpha-i\beta+i\Omega-\lambda-2a^{2}+\frac{C}{2}e^{-(\lambda-i\Omega)\tau}\left(e^{-i(k_{1}-q_{1})}+e^{-i(k_{2}-q_{2})}\right)\right]b_{2} & =0.
\end{align*}
This system has a nontrivial solution if the determinant of
its matrix is zero. We arrive at the characteristic equation
\begin{align} 
%\begin{aligned} 
\chi(\lambda,q_{-},q_{+})= &\; \lambda^{2}+2(a^{2}+R\cos k_{\tau})\lambda+R^{2}+2Ra^{2}\cos k_{\tau}+R_{+}R_{-}e^{-2\lambda\tau+2iq_{+}} - \nonumber  \\
 & -\Bigl[(a^{2}+R\cos k_{\tau}+\lambda)\left(R_{+}e^{ik_{\tau}}+R_{-}e^{-ik_{\tau}}\right) - \label{eq: SLpsolCE} \\
 &\quad\,\, -iR\sin k_{\tau}\left(R_{+}e^{ik_{\tau}}-R_{-}e^{-ik_{\tau}}\right)\Bigr]e^{-\lambda\tau+iq_{+}}, \nonumber
%\end{aligned}
\end{align}
where we additionally denote $q_{\pm}:=\frac{1}{2}\left(q_{1}\pm q_{2}\right)$
and $R_{\pm}:=C\cos(k_{-}\pm q_{-})$. Now the obtained characteristic
equation (\ref{eq: SLpsolCE}) determines the stability of a plane
wave. Namely, for any plane wave, which is defined by the amplitude
$a$, frequency $\Omega$, wave-vectors $k_{+}$ and $k_{-}$ (then
also $k_{\tau}=k_{+}-\Omega\tau$ is given), the equation (\ref{eq: SLpsolCE})
determines the stability with respect to the perturbation mode with
the spatial perturbations $q_{+}$ and $q_{-}$. In particular, if
for all $q_{+}$, $q_{-}\in[0,2\pi]$, all the solutions $\lambda$
of the characteristic equation (\ref{eq: SLpsolCE}) have negative
real parts, then the plane wave is asymptotically stable. The symmetry-relation
$\chi^{\ast}(\lambda^{*},q_{-},-q_{+};-k_{\tau})=\chi(\lambda,q_{-},q_{+};k_{\tau})$
implies, that the stability properties of periodic solutions are invariant
with respect to changing $k_{\tau}\to-k_{\tau}$. Notice that the
obtained equation is a quasi-polynomial, which has generically infinitely many
roots.

Although Eq. (\ref{eq: SLpsolCE}) can be studied numerically for
each given set of parameters, an additional useful analytical insight
in the properties of its solutions is possible under the assumption
of large delay $\tau$. This is performed in the next section.

\subsubsection{Stability of plane wave solutions: large delay approximation}

\subsubsection*{Strong spectrum}

As it was discussed in Sec.~\ref{sub:ststLong}, the strong spectrum
involves only the instantaneous part of the dynamics. Therefore it
does not depend on the spatial perturbation modes $\mathbf{q}$ or
the network-size, since all spatial effects induced by the coupling-structure
are contained in the delayed terms. The reduced characteristic equation
for the strong spectrum can be formally obtained by setting $e^{-\lambda\tau}e^{iq_{+}}=0$
in (\ref{eq: SLststChEq-1}): 
\[
0=\lambda^{2}+2(2a^{2}-\alpha)\lambda+2a^{2}(a^{2}-\alpha)+R^{2}
\]
Its solutions are 
\begin{equation*}
\lambda_{\pm}=\alpha-2a^{2}\pm\sqrt{a^{4}+(a^{2}-\alpha)^{2}-R^{2}}.%\label{eq: SLpsolStrongEV}
\end{equation*}
Any of the solutions $\lambda_{\pm}$ with positive real part belongs
to the strong spectrum. Simple calculations show that there exists
at least one strong eigenvalue with positive real part if 
\begin{equation}
a^{2}<a_{S}^{2}(\alpha;R)=\begin{cases}
\alpha/2 & \mbox{for }-|R|\leq\alpha\leq\sqrt{2}|R|\\
\frac{1}{2}\left(\alpha+\sqrt{\alpha^{2}-2R^{2}}\right) & \mbox{for }\sqrt{2}|R|<\alpha,
\end{cases}
\end{equation}
i.e. if the amplitude of the plane wave is smaller than the one
determined by the curve $a_{S}(\alpha,R)$. Moreover, when the inequality
$a^{4}+(a^{2}-\alpha)^{2}-R^{2}<0$ is satisfied, there are two complex
conjugate unstable eigenvalues $\lambda_{+}=\lambda_{-}^{\ast}\in\mathbb{C}$
with $\operatorname{Re} \lambda_{\pm}=\alpha-2a^{2}$. The bifurcation
diagram in Fig.~\ref{fig:SLBD} illustrates the regions of strong
instability of plane waves (dark gray regions, labeled with $S$).

\subsubsection*{Pseudo-continuous (weak) spectrum}

\begin{figure}
\begin{centering}
\includegraphics[width=0.75\textwidth]{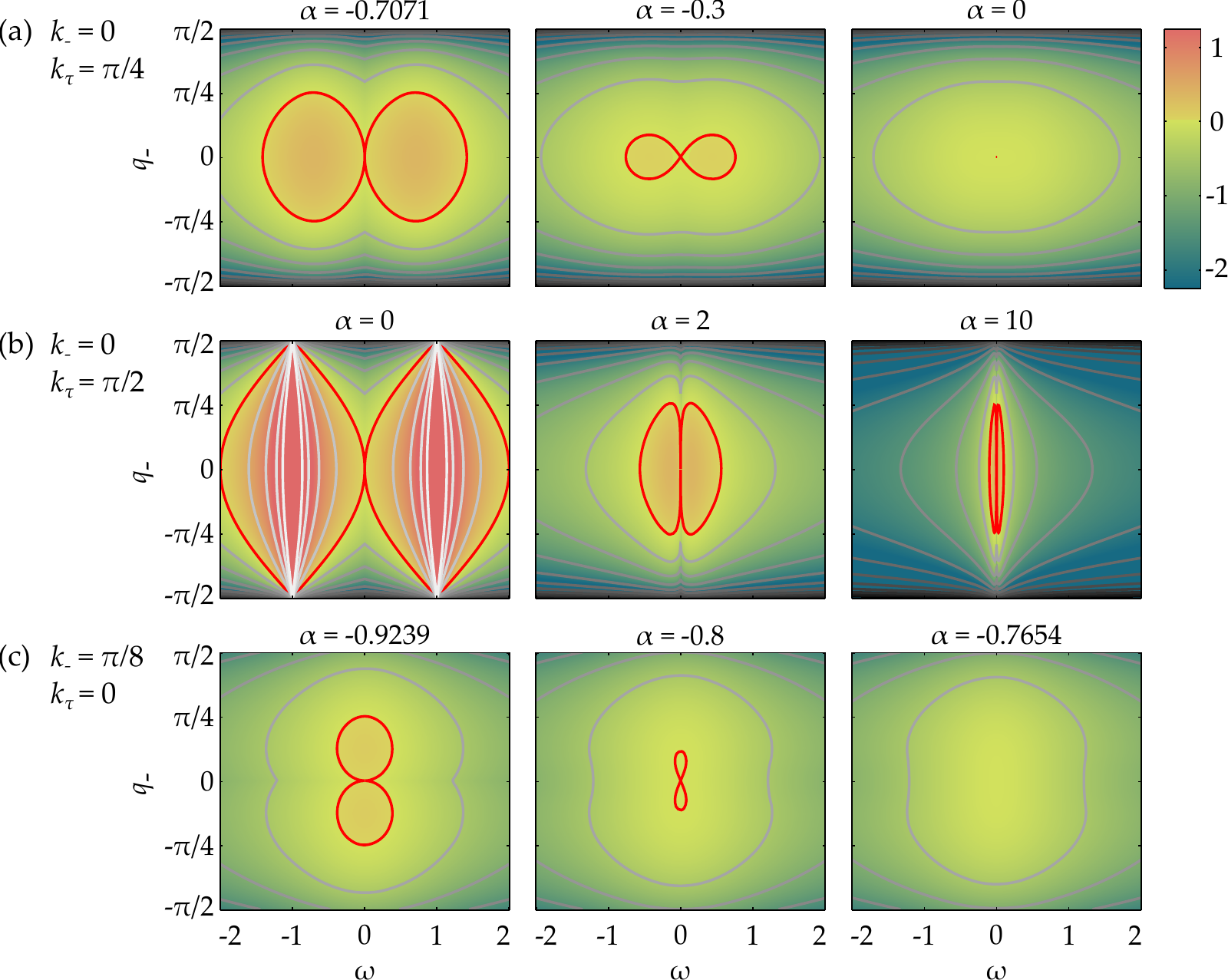}
\par\end{centering}

\caption{Asymptotic continuous spectrum $\gamma_{\mathrm{max}}(\omega,q_{-})$
of periodic solutions for different values of $k_{-}$ and $k_{\tau}$. The density plot shows the scaled real part of the leading eigenvalues. Note that this plot is similar to the Master stability function \cite{Flunkert2010b}.
\label{pic: SL psol pcs}}
\end{figure}

Besides the strong spectrum, there are infinitely many eigenvalues
in the \emph{weak} or \emph{pseudo-continuous spectrum}. Similarly
to Sec.~\ref{sub:ststLong}, this spectrum can be found by substituting
the ansatz 
\[
\lambda=\frac{\gamma(\omega)}{\tau}+i\omega
\]
 into the characteristic equation (\ref{eq: SLpsolCE}). In the limit
of large delay, the terms of the order $\mathcal{O}(1/\tau)$ can
be neglected, resulting in the following equation
\begin{equation}
0=S(q_{-})Y^{2}-2\left[A(\omega,q_{-})+iB(\omega,q_{-})\right]Y+D(\omega)+iE(\omega),\label{eq:psolY}
\end{equation}
with $Y:=e^{-\lambda\tau}e^{iq_{+}}$, and the real valued functions
\begin{align*}
S(q_{-}) & =R_{+}R_{-}=C^{2}\cos{(k_{-}+q_{-})}\cos{(k_{-}-q_{-})},\\
A(\omega,q_{-}) & =C\left(\left[R+a^{2}\cos k_{\tau}\right]\cos{k_{-}}\cos{q_{-}}+\omega\sin{k_{\tau}}\sin{k_{-}}\sin{q_{-}}\right),\\
B(\omega,q_{-}) & =C\left(-a^{2}\sin{k_{\tau}}\sin{k_{-}}\sin{q_{-}}+\omega\cos{k_{\tau}}\cos{k_{-}}\cos{q_{-}}\right),\\
D(\omega) & =R^{2}-\omega^{2}+2Ra^{2}\cos{k_{\tau}},\\
E(\omega) & =2\omega(a^{2}+R\cos{k_{\tau}}).
\end{align*}
Note that the dependence on $q_{+}$ in the obtained equation (\ref{eq:psolY})
is included only in $Y$. Solving the quadratic equation (\ref{eq:psolY})
with respect to $Y$ leads to 
\[
Y_{\pm}(\omega,q_{-})=\frac{1}{S(q_{-})}\left(A(\omega,q_{-})+iB(\omega,q_{-})\pm\sqrt{\zeta(\omega,q_{-})}\right),
\]
with 
\[
\zeta(\omega,q_{-})=A^2(\omega,q_{-})-B^{2}(\omega,q_{-})-S(q_-)D(\omega)+i\left[2A(\omega,q_{-})B(\omega,q_{-})-S(q_-)E(\omega)\right].
\]
Since there are two solutions $Y_{\pm}$, one obtains two functions
for the pseudo-continuous spectrum 
\[
\gamma_{\pm}(\omega,q_{-})=-\log|Y_{\pm}|=-\frac{1}{2}\log\left(Y_{\pm}Y_{\pm}^{\ast}\right).
\]
The spectrum possesses the following symmetries
\begin{equation}
Y_{\pm}(\omega,q_{-}+\pi)=-Y_{\mp}(\omega,q_{-})\label{eq:Ysym1}
\end{equation}
and
\begin{equation}
Y_{\pm}(-\omega,-q_{-})=Y_{\pm}^{\ast}(\omega,q_{-}).\label{eq:Ysym2}
\end{equation}
The first relation (\ref{eq:Ysym1}) implies that it is sufficient
to consider only one of the two functions $\gamma_{\pm}(\omega,q_{-})$,
since they are related to each other by the shift $q_{-}\to q_{-}+\pi$
as 
\begin{equation}
\gamma_{+}(\omega,q_{-}+\pi)=\gamma_{-}(\omega,q_{-}).\label{eq: SL 2-torus gamma shift property}
\end{equation}
This also indicates that the pseudo-continuous spectrum is twofold
degenerate in the limit of $M,N\to\infty$.  The second property
(\ref{eq:Ysym2}) implies that the spectrum has the reflection-symmetry
in the $(\omega,q_{-})$-plane:
\[
\gamma_{\pm}(-\omega,-q_{-})=\gamma_{\pm}(\omega,q_{-}).
\]
Note that in the special case $k_{-}=0$ the additional symmetry-relations
$Y_{\pm}(-\omega,q_{-})=Y_{\pm}^{\ast}(\omega,q_{-})$ and $Y_{\pm}(\omega,-q_{-})=Y_{\pm}(\omega,q_{-})$
hold. Figure~\ref{pic: SL psol pcs} illustrates the function $\gamma_{\max}(\omega,q_{-})=\max\left\{ \gamma_{+},\gamma_{-}\right\} $
for different values of $k_{-}$ and $k_{\tau}$.

The eigenvalues $\lambda$ are known as \emph{characteristic exponents}
or \emph{Floquet-exponents} and are related to the \emph{Floquet-multipliers}
via $\mu=e^{\lambda T}=\exp(\frac{2\pi\gamma}{\Omega\tau})\exp(\frac{2\pi i\omega}{\Omega})$.
As known from the Floquet-theory for periodic solutions, there is
always one trivial multiplier $\mu=1$ or trivial exponent $\lambda=0$,
arising from the continuous symmetry with respect to timeshifts in
autonomous systems (phase shift on the limit cycle). For a perturbation
with $\omega=0$ and $q_{-}=0$ one obtains 
\[
Y_{\pm}(\omega=0,q_{-}=0)=1+\frac{a^{2}}{R}\cos{k_{\tau}}\left(1\pm\frac{\cos{k_{\tau}}}{|\cos{k_{\tau}}|}\right).
\]
The corresponding characteristic exponent follows as 
\begin{align}
\begin{aligned}\gamma_{+}(\omega=0,q_{-}=0)=\begin{cases}
-\log{\left(1+2\frac{a^{2}}{R}\cos{k_{\tau}}\right)}<0 & \mathrm{for}\,\,\cos{k_{\tau}\geq0}\\
0 & \mathrm{for}\,\,\cos{k_{\tau}<0}
\end{cases}\\
\gamma_{-}(\omega=0,q_{-}=0)=\begin{cases}
0 & \mathrm{for}\,\,\cos{k_{\tau}\geq0}\\
-\log{\left(1+2\frac{a^{2}}{R}\cos{k_{\tau}}\right)}>0 & \mathrm{for}\,\,\cos{k_{\tau}<0}
\end{cases}
\end{aligned}
\label{eq:SL2-torustrivialmultiplier}
\end{align}
Note that this property of the spectrum is not affected by the long
delay approximation, since the approximation becomes exact at $\gamma=0$.
Apparently there are two parameter domains separated by $\cos{k_{\tau}}=0$.
Using (\ref{eq: SLpsolRe}), one finds that this boundary corresponds
to the curve $a_{U}:=\sqrt{\alpha}$. Thus, all periodic solutions
with the amplitudes smaller than $a_{U}$ for a given $\alpha$ are
unstable due to a positive characteristic exponent with $\omega=0$.
According to \cite{Cross1993} this instability is called a \emph{uniform
instability}. In order to determine the neutral stability curve, the
following discussion is restricted to the regime with $\cos{k_{\tau}}\geq0$.
Since the relation (\ref{eq: SL 2-torus gamma shift property}) holds,
we will focus the analysis on $\gamma_{-}(\omega,q_{-})$. 

The trivial multiplier always denotes a critical point of the pseudo-continuous
spectrum at $(\omega=0,q_{-}=0)$, where the gradient vanishes: 
\[
\left.\nabla\gamma_{-}\right|_{\omega=0,q_{-}=0}=\left.\left(\frac{\partial}{\partial\omega},\frac{\partial}{\partial q_{-}}\right)\gamma_{-}\right|_{\omega=0,q_{-}=0}=0.
\]
This can be verified by a direct calculation. Therefore the point
$(\omega=0,q_{-}=0)$ is either an extremum or saddle of the pseudo-continuous
spectrum. Analyzing the shape of the spectrum close to the trivial
multiplier shows the appearance of the \emph{modulational instability}
\cite{Cross2009,Cross1993}, see also Fig.~\ref{pic: SL psol pcs}.
For this, let us consider the second order approximation of $\gamma_{-}$
at $(\omega,q_{-})=0$, involving the Hessian matrix 
\[
H=\left.\left(\begin{array}{cc}
\frac{\partial^{2}\gamma_{-}}{\partial\omega^{2}} & \frac{\partial^{2}\gamma_{-}}{\partial\omega\partial q_{-}}\\
\frac{\partial^{2}\gamma_{-}}{\partial\omega\partial q_{-}} & \frac{\partial^{2}\gamma_{-}}{\partial q_{-}^{2}}
\end{array}\right)\right|_{\omega=0,q_{-}=0}.
\]
Direct calculation leads to the following expressions for the elements
of the Hessian matrix 
\begin{align*}
\left.\frac{\partial^{2}\gamma_{-}}{\partial\omega^{2}}\right|_{\omega=0,q_{-}=0} & =\frac{1}{R^{2}\cos^{2}{k_{\tau}}}\left(\frac{R}{a^{2}}\frac{\sin^{2}{k_{\tau}}}{\cos{k_{\tau}}}-1\right),\\
\left.\frac{\partial^{2}\gamma_{-}}{\partial q_{-}^{2}}\right|_{\omega=0,q_{-}=0} & =-1+\frac{R\tan^{2}{k_{-}}}{a^{2}\cos^{3}{k_{\tau}}}+\tan^{2}{k_{-}}\tan^{2}{k_{\tau}},\\
\left.\frac{\partial^{2}\gamma_{-}}{\partial\omega\partial q_{-}}\right|_{\omega=0,q_{-}=0} & =\frac{\tan{k_{-}}\tan{k_{\tau}}}{a^{2}\cos^{2}{k_{\tau}}}.
\end{align*}
The curvature of the asymptotic continuous spectrum close to the trivial multiplier is directly
related to the stability of the corresponding plane wave. If the surface
is locally concave close to $(\omega,q_{-})=(0,0)$, then the Hessian
is negative definite and the corresponding periodic orbit is stable
(at least the part of the spectrum which is close to $(\omega,q_{-})=(0,0)$).
Otherwise, if the curvature is convex (Hessian is positive definite)
or the origin is a saddle-point (Hessian is indefinite), the plane
wave is unstable. The curvature is characterized by the real eigenvalues
of the symmetric Hessian matrix. The analysis of the eigenvalues of
the Hessian matrix leads to the following condition for the stability
of the plane wave, which is the condition for the negativeness of
the eigenvalues of $H$:
\[
\left(\cos^{2}{k_{\tau}}-\sin^{2}{k_{-}}\right)\left(R\cos{k_{\tau}}+a^{2}\right)-R\cos{k_{\tau}}>0.
\]
Using the amplitude relation (\ref{eq: SLpsolRe}), the bifurcation
is described by a 3rd order polynomial in $a^{2}$: 
\begin{equation}
0=a^{6}-\frac{5}{2}\alpha a^{4}+\left(2\alpha^{2}-\frac{R^{2}}{2}(1+2\sin^{2}{k_{-}}\right)a^{2}-\frac{\alpha^{3}}{2}+\frac{R^{2}\alpha}{2}(1+\sin^{2}{k_{-}}).\label{eq: SL 2-torus general neutral stability curve}
\end{equation}
Solving the corresponding to (\ref{eq: SL 2-torus general neutral stability curve})
equation for $a^{2}$ gives the neutral stability curve for an arbitrary
plane wave with $|k_{-}|\in[0,\frac{\pi}{2}]$. The analytical solution
can be found by using Cardano's method, but is not written here for
brevity. In the particular case $k_{-}=0$, the neutral stability
curve can be simply expressed as 
\begin{equation}
a_{M}^{2}=\frac{1}{4}\left(3\alpha+\sqrt{\alpha^{2}+8C^{2}}\right),\label{eq: SL 2-torus Eckhaus curve}
\end{equation}
which coincides with the result obtained in \cite{Yanchuk2008a}
for the ring of coupled oscillators without delay (and also for a ring
with large delay). By substituting (\ref{eq: SLpsolRe}) into (\ref{eq: SL 2-torus general neutral stability curve}),
one can obtain the minimal $\alpha=\alpha_{0}$ with
\[
\alpha_{0}(k_{-},k_{\tau})=R\cos k_{\tau}\frac{1-2(\cos^{2}k_{-}-\sin^{2}k_{\tau})}{\cos^{2}k_{-}-\sin^{2}k_{\tau}},
\]
where a plane wave with particular $k_{-}$ and $k_{\tau}$ stabilizes.

\begin{figure}
\centering \includegraphics[width=0.75\linewidth]{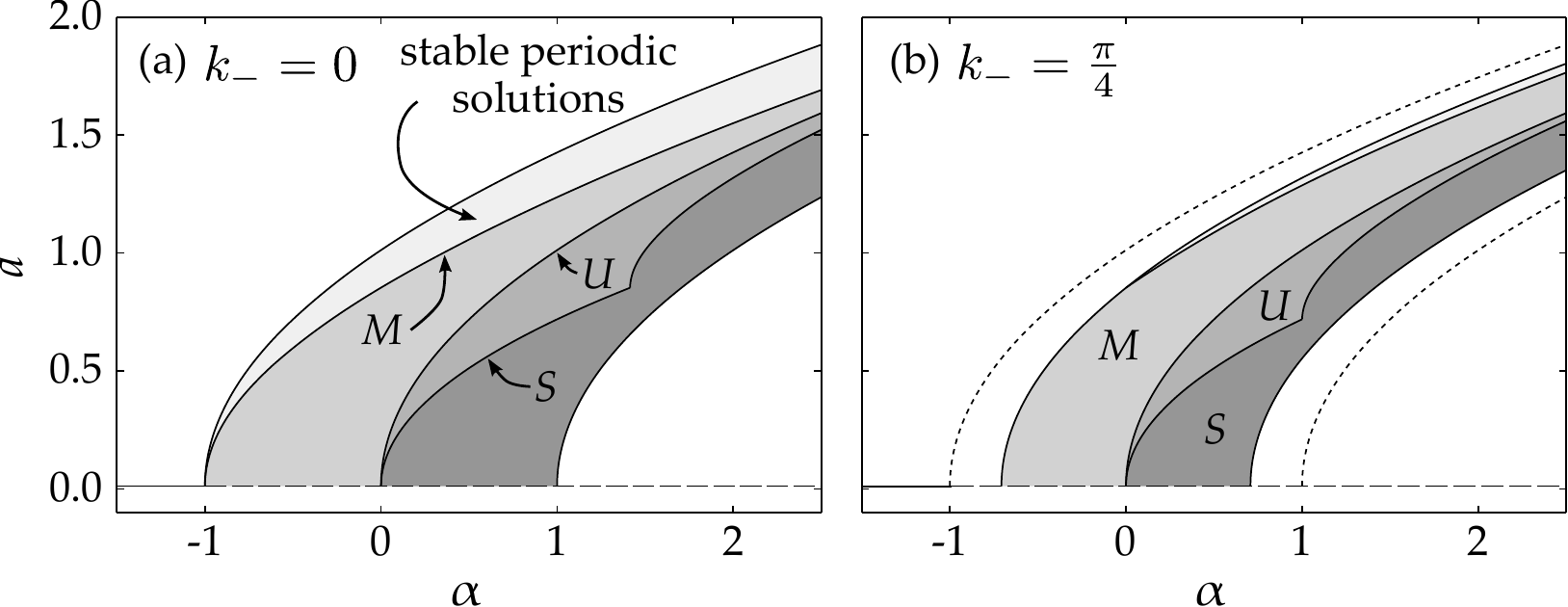}
\caption{Stability diagrams for periodic solutions of the SL-system (\ref{eq: SL})
with different values of $k_{-}$: (a) $k_{-}=0$ and (b) $k_{-}=\pi/4$.
The regime of stable periodic solutions is filled with a light gray
and is located above the $M$-labeled curve in this projection. The
curves denote the upper boundaries of the different stability regimes
discussed in the text: $M$ modulational instability, $U$ uniform
instability and $S$ strong instability. The size of the stable regime
decreases with the increasing of $k_{-}$. \label{fig:SLBD}}
\end{figure}

The bifurcation diagram in Fig.~\ref{fig:SLBD} summarizes and illustrates
the obtained results, showing the regions where the plane waves are
stable (light gray), weakly unstable (darker gray, labeled with $U$
and $M$), and strongly unstable (dark gray, labeled with $S$).

\subsubsection{Qualitative summary for the plane waves in coupled Stuart-Landau
model}

The main qualitative conclusions of this relatively technical section
are as follows: The family of plane wave solutions (\ref{eq:persol})
is located on the circles (\ref{eq: SLpsolCir}) (for a fixed $k_{-}$),
and the number of plane waves grows accordingly to (\ref{eq:Npsol}).
The stability of a plane wave is governed by the characteristic equation
(\ref{eq: SLpsolCE}) and determined by its position on the circle.
More specifically, the plane waves with the higher amplitude tend
to be more stable than those with the lower amplitude. Figure \ref{fig:SLBD}
illustrates this by showing stable, as well as weakly and strongly
unstable ``positions'' on the circle. Thus, with the increasing
$\alpha$, the number of stable plane waves increases. Plane waves
with smaller $\left|k_{-}\right|$ also tend to be more stable
than those with larger $\left|k_{-}\right|$. Therefore we expect
that the plane waves which are almost diagonal are more abundantly
observed.

\section{FitzHugh-Nagumo neurons with homogeneous coupling delays\label{sec:FHN-hom}}

In this section we consider a lattice of $M\times N$ delay-coupled
FitzHugh-Nagumo neurons, which are coupled via excitatory chemical
synapses. The coupling architecture is the same as described in section
\ref{sec: model system}. The model system reads 
\begin{align}  \label{eq:FHN}
\begin{aligned}
\frac{\mathrm{d}}{\mathrm{d}t}v_{m,n} &= v_{m,n}-\frac{1}{3}v_{m,n}^{3}-w_{m,n}+I+\\
&\quad +\frac{C}{2}(v_{r}-v_{m,n})(s_{m-1,n}(t-\tau)+s_{m,n-1}(t-\tau))\\
\frac{\mathrm{d}}{\mathrm{d}t}w_{m,n} &= \varepsilon(v_{m,n}+a-bw_{m,n})\\
\frac{\mathrm{d}}{\mathrm{d}t}s_{m,n} &= \alpha(v_{m,n})(1-s_{m,n})-0.6s_{m,n}
\end{aligned}
\end{align}
with $\alpha(v)=\frac{1}{2}[1+e^{-5(v-1)}]^{-1}$. The variable $v_{m,n}$
denotes the membrane potential of the corresponding neuron and $w_{m,n}$
is a recovery variable, combining several microscopic dynamical variables
of the biological neuron. The external stimulus current applied to
the neuron is denoted by $I$ and $C$ is the coupling strength. We
fix the parameters $a=0.7$, $b=0.8$, and $\varepsilon=0.08$. The
reversal potential is taken as $v_{r}=2$ for excitatory coupling.
Similar model equations have been investigated in \cite{Perlikowski2010,Yanchuk2011}
for unidirectional rings.

We demonstrate that the destabilization of the homogeneous steady
state, the set of plane waves as well as their stability properties
possess the same qualitative features which we observed in the previous
section on Stuart-Landau oscillators. However, the apparent difficulty
for the analysis of nonlinear plane waves is that they are not known
analytically. Therefore we use numerical bifurcation analysis with
\textsc{DDE-Biftool} \cite{Engelborghs2002} and have to restrict ourselves to
relatively small lattice size and delay values in Sec.~\ref{sub:FHNH}.

\subsection{Homogeneous steady state and its stability}

\begin{figure}[t]
\centering \includegraphics[width=1\linewidth]{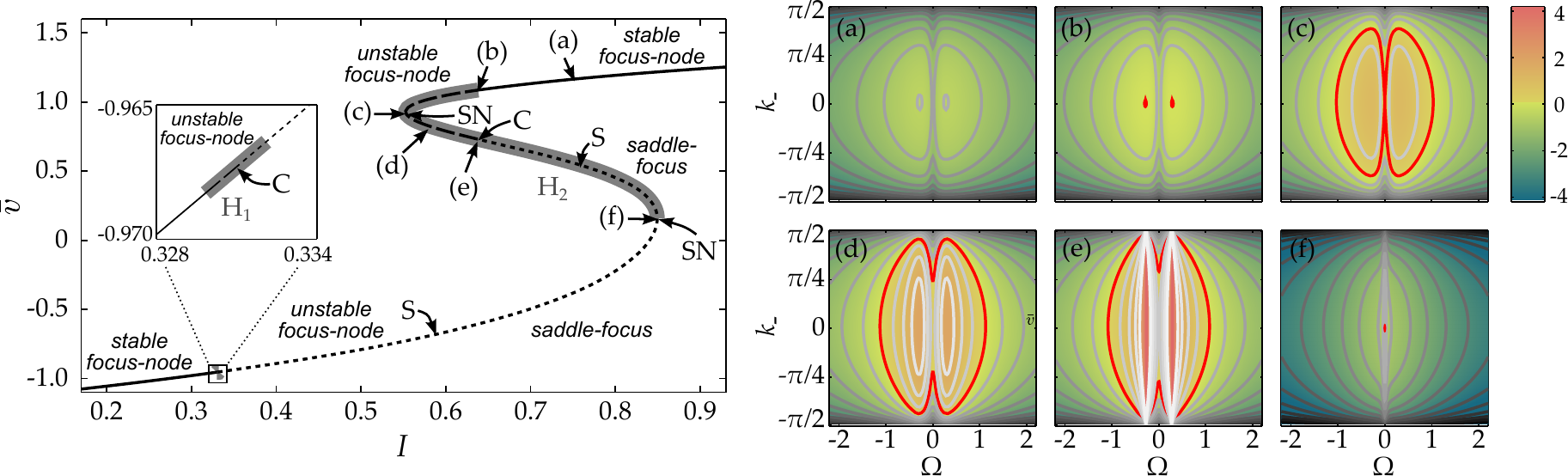}
\caption{Homogeneous stationary state and its bifurcations for $C=3$ for the FitzHugh-Nagumo model Eq. (\ref{eq:FHN}). In the
case of an infinitely large lattice, the spectrum is described by
the hybrid dispersion relation (\ref{eq: FHNststhybrid}), parametrized
by two continuous parameters $\Omega$ and $k_{-}$, shown on the
right. The color denotes the real parts of the eigenvalues: (a) stable
stationary state, (b) destabilization via Hopf-bifurcation, (c) saddle-node
bifurcation (SN), where an eigenvalue with $\Omega=0$ becomes unstable,
(d) weakly unstable stationary state, (e) cusp bifurcation (C) and
appearance of a strongly unstable spectrum, (f) saddle-node bifurcation
and lower boundary of the Hopf-domain $H_{2}.$ \label{fig:FHMststhyb} }
\end{figure}

The system (\ref{eq:FHN}) has a homogeneous steady state $\bar{\mathbf{u}}=\left(\bar{v},\bar{w},\bar{s}\right)$.
The value for the membrane resting potential $\bar{v}$ can be obtained
as a solution of the scalar equation 
\begin{equation}
0=\bar{v}-\frac{1}{3}\bar{v}^{3}-\frac{\bar{v}+a}{b}+I+C(v_{r}-\bar{v})\frac{\alpha(v)}{\alpha(v)+0.6}.\label{eq:FHNstst}
\end{equation}
The steady-state values of the remaining variables follow as 
\[
\bar{w}=\frac{\bar{v}+a}{b},\qquad\mathrm{and}\qquad\bar{s}=\frac{\alpha(v)}{\alpha(v)+0.6}.
\]
In the case of weak coupling strength the homogeneous stationary state
is unique, but for $C_{\mathrm{SN}}=1.46475$ a saddle-node bifurcation
of the equilibrium takes place. For strong coupling $C>C_{\mathrm{SN}}$
there is a domain of the control parameter $I$ with three coexisting
stationary states, see Fig. \ref{fig:FHMststhyb}. In order to analyze
the stability of the stationary state, we derive the linearized evolution
equation for small perturbations of the equilibrium and subsequently
diagonalize it in Fourier-space, just as in the previous section.
One obtains the system 
\[
\frac{\mathrm{d}}{\mathrm{d}t}\mathbf{\delta\tilde{\mathbf{u}}}(t)=A\mathbf{\delta\tilde{u}}(t)+2B\cos{\left(k_{-}\right)}e^{ik_{+}}\mathbf{\delta\tilde{u}}(t-\tau)
\]
with the real valued matrices 
\[
A=\left(\begin{array}{ccc}
1-\bar{v}^{2}-C\bar{s} & -1 & 0\\
\varepsilon & -b\varepsilon & 0\\
5\alpha(\bar{v})(1-2\alpha(\bar{v}))(1-\bar{s}) & 0 & -\alpha(\bar{v})-0.6
\end{array}\right)
\]
and 
\[
B=\left(\begin{array}{ccc}
0 & 0 & \frac{C}{2}(v_{r}-\bar{v})\\
0 & 0 & 0\\
0 & 0 & 0
\end{array}\right).
\]
The corresponding characteristic equation reads 
\begin{equation}
0=\prod_{\{(k_{+},k_{-})\}}\det{\left|-\lambda\operatorname{Id}+A+2B\cos{\left(k_{-}\right)}e^{ik_{+}}e^{-\lambda\tau}\right|},\label{eq:FHNcheq}
\end{equation}
where $\operatorname{Id}$ represents the $3\times3$ identity matrix.
Similarly to the previous analysis, the stability of the homogeneous
steady state is completely determined by the Eq.~(\ref{eq:FHNcheq}),
which can be studied numerically using e.g. Newton-Raphson iteration.
An additional insight in the properties of the spectrum can be given
using the large delay approximation, which is done in Sec~\ref{sub:FHNlong}.
Further, in Sec.~\ref{sub:FHNH}, we study the appearance of periodic
solutions in Hopf bifurcations.

\subsubsection{Large delay approximation} \label{sub:FHNlong}

\subsubsection*{Strongly unstable spectrum}

The strongly unstable spectrum results from considering only the instantaneous
part of Eq. (\ref{eq:FHNcheq}) 
\[
0=\det{\left|A-\lambda\operatorname{Id}\right|}=(a_{33}-\lambda)\left[(a_{11}-\lambda)(a_{22}-\lambda)-a_{21}a_{12}\right].
\]
There is always one real solution $\lambda_{0}=a_{33}\leq-0.6$ of
this 3rd-order polynomial, which is strictly negative. The remaining
eigenvalues are 
\[
\lambda_{\pm}=\frac{1}{2}\left(a_{11}-b\varepsilon\pm\sqrt{(a_{11}+b\varepsilon)^{2}-4\varepsilon}\right).
\]
Note that all relevant parameters are contained in $a_{11}=1-\bar{v}^{2}-C\bar{s}$,
involving the current $I$ and coupling strength $C$ directly or
via the corresponding homogeneous steady state, respectively. 

The strongly unstable spectrum exists when the real part of the largest
eigenvalue $\lambda_{+}$ is positive. This is the case when $a_{11}>b\varepsilon$.
The appearance of the strongly unstable spectrum occurs at the cusp-bifurcation
of the asymptotic continuous spectrum and labled as ``C'' in Fig. \ref{fig:FHMststhyb}. Moreover,
there exists a pair of complex conjugate eigenvalues for $a_{11}<2\sqrt{\varepsilon}-b\varepsilon$,
otherwise the eigenvalues are real. The corresponding boundary is
labeled with ``S'' in Fig. \ref{fig:FHMststhyb} and mediates the
transition between an unstable focus-node and a saddle-focus.

\subsubsection*{Pseudo-continuous spectrum}

The primary bifurcations of the steady state are captured by the pseudo-continuous
spectrum. Just as in the previous section in the case of Stuart-Landau
oscillators, this can be found by applying the ansatz $\lambda=\gamma/\tau+i\Omega$.
Neglecting terms of order $\mathcal{O}(1/\tau)$ and introducing the
variable $Y=e^{-\gamma}e^{i(k_{+}-\Omega\tau)}$ yields the modified
characteristic equation 
\begin{align*}
0 & =\det{\left|-i\Omega\operatorname{Id}+A+2B\cos k_{-}Y\right|}=\det{\left|\begin{array}{ccc}
a_{11}-i\Omega & a_{12} & 2b_{13}\cos{k_{-}}Y\\
a_{21} & a_{22}-i\Omega & 0\\
a_{31} & 0 & a_{33}-i\Omega
\end{array}\right|} = \\
 & =(a_{11}-i\Omega)(a_{22}-i\text{\ensuremath{\Omega}})(a_{33}-i\Omega)-2(a_{22}-i\Omega)a_{31}b_{13}\cos{k_{-}}Y-(a_{33}-i\Omega)a_{12}a_{21}.
\end{align*}
Due to the simple linear coupling-structure, this is a linear equation
in $Y$, which can be solved as 
\[
Y=\frac{a_{33}-i\Omega}{2a_{31}b_{13}\cos{k_{-}}}\left(a_{11}-i\Omega-\frac{a_{12}a_{21}}{a_{22}-i\Omega}\right),
\]
leading to the asymptotic continuous spectrum 
%\begin{align} \label{eq: FHNststhybrid}
%\begin{aligned}
%\gamma(\Omega,k_{-}) & =-\log{|Y(\Omega,k_{-})|}%\\
% & =-\frac{1}{2}\log{\left[\frac{a_{33}^{2}+\Omega^{2}}{(2b_{13}a_{31}\cos{k_{-}})^{2}}\left(a_{11}^{2}+\Omega^{2}+\frac{a_{12}a_{21}}{a_{22}^{2}+\Omega^{2}}\left(a_{12}a_{21}-2a_{11}a_{22}+2\Omega^{2}\right)\right)\right]}.
%\end{aligned}
%\end{align}
\begin{equation} \label{eq: FHNststhybrid}
\gamma(\Omega,k_{-}) = -\log{|Y(\Omega,k_{-})|}
\end{equation}

This is a function of two parameters, determining the spectrum and
stability of the steady state with respect to the perturbations with
the spatial mode $k_{-}$ (and independent on $k_{+})$ and the delay-induced temporal modes $\Omega$. Some plots of this surface are illustrated
in Fig.~\ref{fig:FHMststhyb}. Apparently the destabilization is
similar to the case of Stuart-Landau oscillators described in Sec. \ref{sub:ststLong}.
The asymptotic weak spectrum is invariant with respect to $\Omega\to-\Omega$,
$k_{-}\to-k_{-}$ and $k_{-}\to k_{-}+n\pi$, $n\in\mathbb{Z}$. The
bifurcations of (\ref{eq: FHNststhybrid}) lead to the boundaries
of domains, where Hopf bifurcations are possible (shown as $H_{1}$
and $H_{2}$ in Fig. \ref{fig:FHMststhyb}), and saddle-node-bifurcations.
Many properties (such as extrema and roots) of the hybrid dispersion
relation (\ref{eq: FHNststhybrid}) are analytically accessible, but not
given here explicitly, since they involve solutions of 3rd order polynomials.

\subsection{Hopf bifurcations and periodic attractors\label{sub:FHNH}}

\begin{figure}[th]
\centering \includegraphics[width=0.5\textwidth]{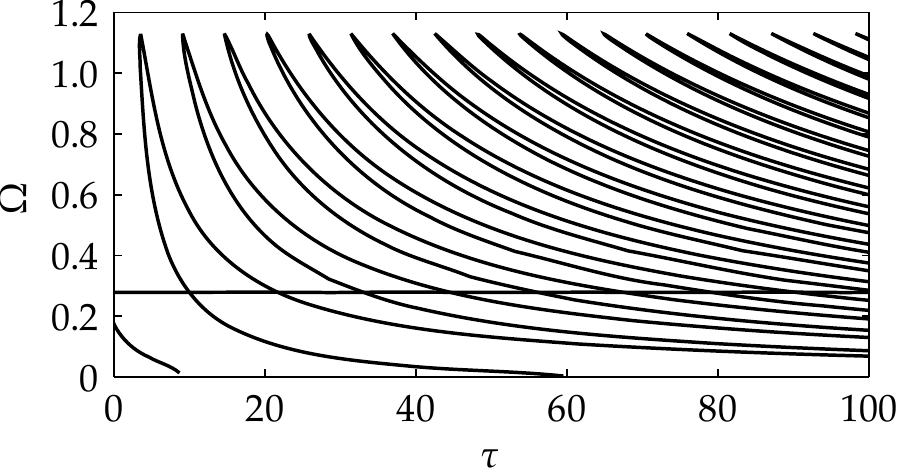}
\caption{Hopf bifurcation branches of a single FHN oscillator with delayed
feedback with $C=3$. For $M\times N$ lattices similar structures
can be obtained. The parameters used are given in the text.\label{pic: FHN Hopf branches}}
\end{figure}

The Hopf bifurcations of the stationary state can be computed by setting
$\lambda=i\Omega$ in Eq. (\ref{eq:FHNcheq}) and subsequently separating
the resulting complex equation into two real equations. Using $\bar{v}$
and $\omega$ as parameters, the system can be solved by using e.g.
Newton-Raphson iterations. The corresponding value of the bifurcation-parameter
$I$ follows uniquely from Eq. (\ref{eq:FHNstst}). Using the software
package \textsc{DDE-Biftool} \cite{Engelborghs2002},
we perform a continuation of the Hopf bifurcations in the $(I,\tau)$-plane.
The result is shown in Fig. \ref{pic: FHN Hopf branches}, where the
Hopf frequency $\Omega$ is plotted vs. the time delay $\tau$. The
structure of the branches can be understood by using reappearance
arguments for periodic solutions \cite{Yanchuk2009}. Some of the
Hopf branches terminate with zero frequency in a homoclinic bifurcation.

\begin{figure}[h]
\centering \includegraphics[width=1\linewidth]{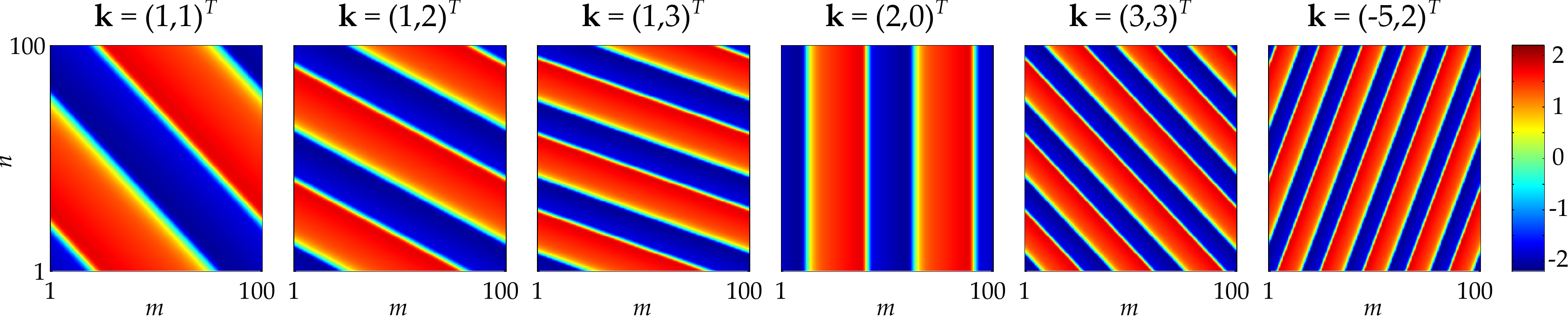}
\caption{Snapshots of several coexisting stable travelling waves in a $100\times100$-lattice
of delay-coupled FHN neurons with $C=\text{3}$, $I=0$ and $\tau=50$.
The color denotes the value of the membrane voltage $v_{m,n}$ of the corresponding
neuron.\label{pic: FHN travelling waves snapshots}}
\end{figure}

Using \textsc{DDE-Biftool}, we perform a numerical continuation of
the periodic solutions, emerging from the Hopf bifurcations. The spatial
orientation of a periodic solution is conserved along the branch,
while varying the external current $I$ as a control parameter. Typically,
a periodic solution connects two Hopf points, which are both solutions
of Eq. (\ref{eq:FHNcheq}) with the same $k_{+}$ and $k_{-}$. For
vanishing delay, all stable periodic orbits are diagonal travelling
waves with $k_{-}=0$, including the synchronized solution. Increasing
the coupling delay significantly enhances the stability properties
of periodic solutions and allows for stable travelling waves with $k_{-}\neq0$.
Moreover, the periodic solutions appear in a larger regime of the
control parameter $I$. Snapshots of several coexisting travelling
waves in a system of $100\times100$ FHN-neurons with $\tau=50$ are
shown in Fig. \ref{pic: FHN travelling waves snapshots}. Such solutions serve as
the starting point for the more complicated patterns in systems with
inhomogeneous delays, discussed in the following section.

\begin{figure}[h]
\centering \includegraphics[width=0.7\linewidth]{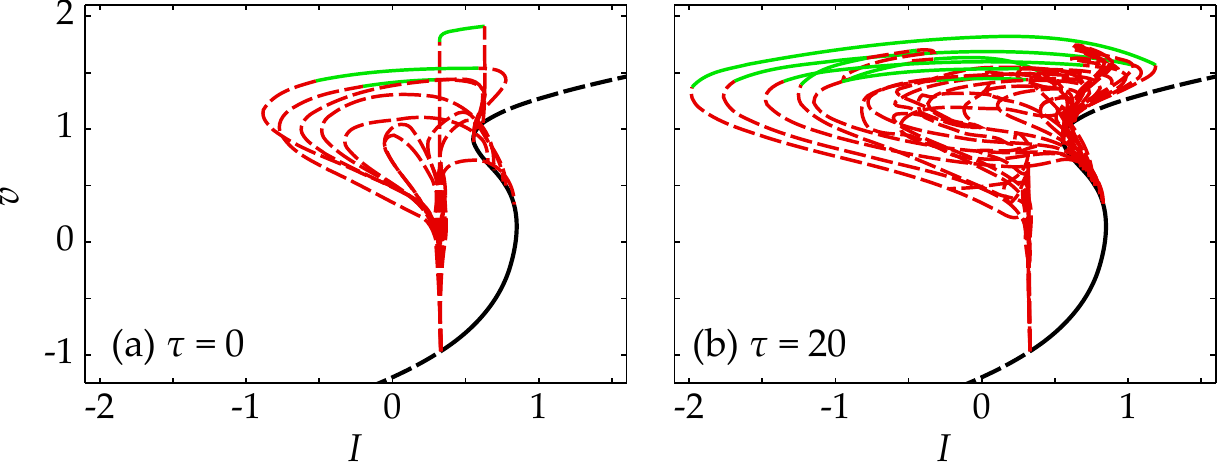}
\caption{Bifurcation diagrams for $3\times3$-lattices of coupled FHN neurons
with $C=\text{3}$ and different delays: (a) $\tau=0$ and (b) $\tau=20$.
Green solid lines denote stable periodic solutions, red dashed lines
show unstable ones. The stable stationary state is depicted by a black
solid line, unstable steady states by a black dashed line. The parameters used are given in the text.\label{pic: FHN bifurcationdiagram}}
\end{figure}

\section{Patterns in systems with inhomogeneous delays\label{sec:non}}

\subsection{Componentwise time-shift transformation}

Consider a delayed dynamical system with a coupling topology as described
by Eq. (\ref{eq: general system}) with homogeneous delays. As follows
from the previous analysis, for sufficiently large delay $\tau$ the
homogeneous system possesses a large set of stable plane waves, including
stable spatially homogeneous time-periodic solutions. 

We can transform the plane waves of the homogeneous system into an
(almost) arbitrary pattern by adjusting the coupling delays. The derivation
of the transformation presented here is a generalization of the method
described in \cite{Popovych2011,Yanchuk2011} for unidirectionally
coupled rings and \cite{Luecken2013b} for arbitrary networks with
delays.

Rewriting the system (\ref{eq: general system}) with homogeneous
delays $\tau$ with respect to the new coordinates $\mathbf{v}$ given
by $\mathbf{u}_{m,n}(t)=\mathbf{v}_{m,n}(t-\eta_{m,n})$ leads to
the system 
%\begin{align*}
%\dot{\mathbf{v}}_{m,n}(t) & =\mathbf{F}\left(\mathbf{v}_{m,n}(t),\,\mathbf{v}_{m-1,n}(t+\eta_{m,n}-\eta_{m-1,n}-\tau)+\mathbf{v}_{m,n-1}(t+\eta_{m,n}-\eta_{m,n-1}-\tau)\right)\\
% & =\mathbf{F}\left(\mathbf{v}_{m,n}(t),\,\mathbf{v}_{m-1,n}(t-\tau_{m,n}^{\downarrow})+\mathbf{v}_{m,n-1}(t-\tau_{m,n}^{\rightarrow})\right)
%\end{align*}
\begin{align*}
\dot{\mathbf{v}}_{m,n}(t) &= \mathbf{F} \big( \mathbf{v}_{m,n}(t),\,\mathbf{v}_{m-1,n}(t+\eta_{m,n}-\eta_{m-1,n}-\tau) + \\
&\qquad \qquad \qquad \, +  \mathbf{v}_{m,n-1}(t+\eta_{m,n}-\eta_{m,n-1}-\tau) \big) \\
&= \mathbf{F}\left(\mathbf{v}_{m,n}(t),\,\mathbf{v}_{m-1,n}(t-\tau_{m,n}^{\downarrow})+\mathbf{v}_{m,n-1}(t-\tau_{m,n}^{\rightarrow})\right)
\end{align*}
with adjusted non-homogeneous delays

\begin{align} \label{eq:delaytransformation}
\begin{aligned}
\tau_{m,n}^{\downarrow} &=  \tau-\eta_{m,n}+\eta_{m-1,n},\\
\tau_{m,n}^{\rightarrow}&=  \tau-\eta_{m,n}+\eta_{m,n-1}.
\end{aligned}
\end{align}
The timeshifts $\eta\in\mathbb{R}^{M\times N}$ can have an arbitrary
form, up to the restriction that the new delays need to be positive.
By adjusting the timeshifts, one can obtain stable, time-periodic
attractors of various spatial forms. For example, a stable synchronous
periodic solution $\mathbf{u}_{m,n}(t)=\mathbf{u}_{0}(t)=\mathbf{u}_{0}(t+T)$
of the homogeneous system corresponds to a solution $\mathbf{v}_{m,n}(t)=\mathbf{u}_{0}(t+\eta_{m,n})$
of the non-homogeneous system, where each component is shifted in
time by an amount $\eta_{m,n}$. According to \cite{Luecken2013a,Luecken2013b},
the stability properties of the periodic solutions are kept by the
transformation. This shift will result to a shifted value of the dynamical
variables (e.g. voltage for the neuronal models). Thus, the encoded
pattern $\eta_{m,n}$ will be visible in the dynamical variables.
Since $\eta_{m,n}$ is practically arbitrary, there is a variety of
patterns, which can appear as stable attractors in the systems with
inhomogeneously delayed connections.

\subsection{Examples\label{sub:Ex}}

\begin{figure}[h]
\centering \includegraphics[width=1\textwidth]{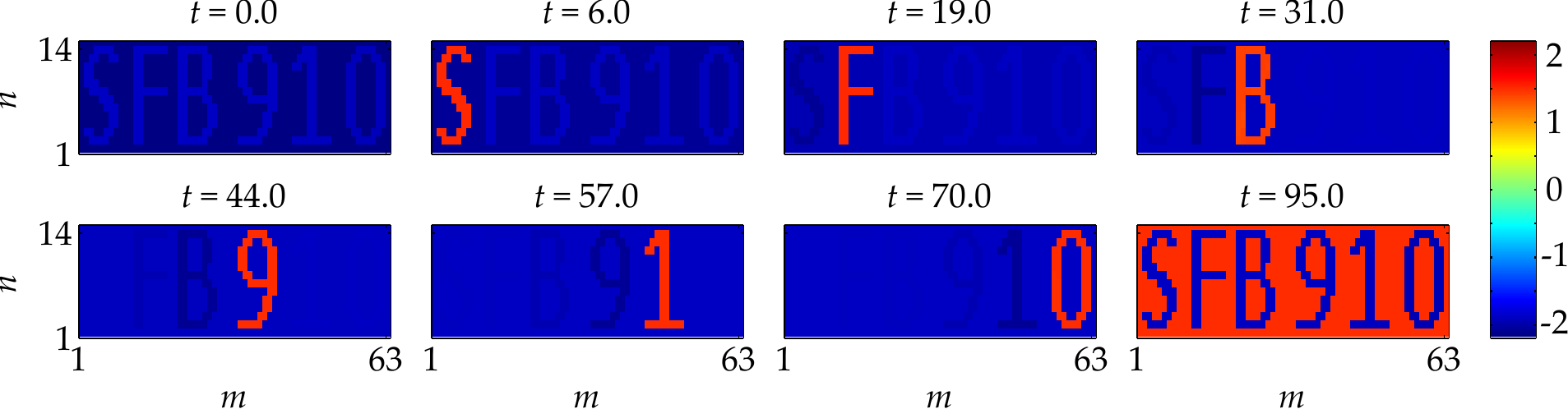}
\caption{Snapshots of the ``SFB 910''-pattern in a $63\times14$-lattice
of delay-coupled FHN-neurons with $C=\text{3}$ and $I=-2$. There
are several clusters of neurons, spiking synchronously. The pattern
oscillates with a period $T=102.3$.\label{pic: FHN nonhom SFB} }
\end{figure}

\begin{figure}[h]
\centering \includegraphics[width=1\textwidth]{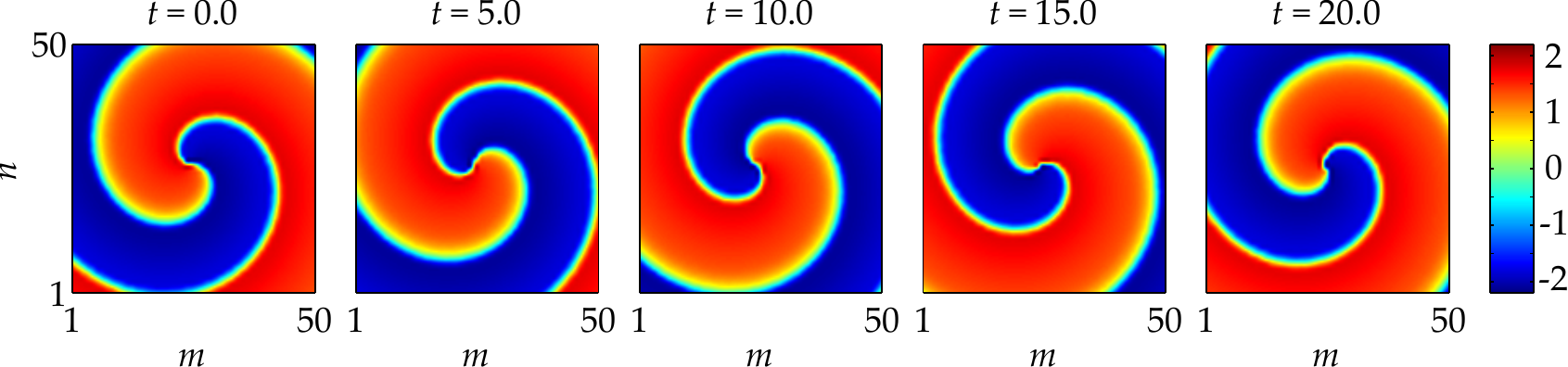}
\caption{Snapshots of a spiral wave pattern in a $50\times50$-lattice of delay-coupled
FHN-neurons with $C=\text{3}$ and $I=0$. The pattern periodically
reappears with a period $T=25.8$.\label{pic: FHN nonhom spiral}}
\end{figure}

Some illustrative examples for stable spatio-temporal patterns in a lattice with non-homogeneous
delays are shown in Figs.~\ref{fig:ML}, \ref{pic: FHN nonhom SFB},
and \ref{pic: FHN nonhom spiral}. All examples are constructed from
synchronized solutions with $\mathbf{k}=(0,0)^{T}$ via the the delay-transformation
(\ref{eq:delaytransformation}). However, the scaling of the patterns
$\eta_{m,n}$ with respect to the period time is different in the
examples. In the ``Mona Lisa''-pattern (Fig. \ref{fig:ML}) the spiking times are chosen
only slightly different, so that the pattern is a slightly adapted
standing front solution. In the examples in Fig. \ref{pic: FHN nonhom SFB}
and \ref{pic: FHN nonhom spiral} the spiking-times are distributed
over the whole period.

\section{Conclusions}
The linear stability of time-periodic patterns has been studied analytically in the limit of large delay in terms of a hybrid dispersion relation.

We have shown that arbitrary stable spatio-temporal periodic patterns can be created in two-dimensional lattices of coupled oscillators with inhomogeneous coupling delays. We propose that this offers 
interesting applications for the generation, storage, and information processing of visual patterns, 
for instance in networks of optoelectronic \cite{Williams2013} or electronic \cite{Rosin2013} oscillators.
Our results have been illustrated with two models of the local node dynamics which have a wide range of applicability: (i) the Stuart-Landau oscillator, i.e., a generic model which arises by center-manifold expansion of a limit cycle system near a supercritical Hopf bifurcation, and (ii) the FitzHugh-Nagumo model, which is a generic model of neuronal spiking dynamics.

\section*{Acknowledgements}

We thank L. L{\"u}cken and M. Zaks for useful discussions and the DFG
for financial support in the framework of the Collaborative Research
Center SFB 910. 

\bibliographystyle{siam}
\bibliography{2dlattice}

\end{document}